\documentclass[12pt]{amsart}
\usepackage[top=1in, bottom=1in, left=1in, right=1in]{geometry}

\usepackage{amssymb,amsmath,amsthm}
\usepackage[hyphens]{url}

\usepackage{graphicx}
\usepackage{color}

\numberwithin{equation}{section}

\newtheorem{mainthm}{Theorem}
\newtheorem{mainconj}{Conjecture}
\newtheorem{mainprop}{Proposition}

\newtheorem{thm}{Theorem}[section]
\newtheorem{lem}[thm]{Lemma}
\newtheorem{prop}[thm]{Proposition}
\newtheorem{cor}[thm]{Corollary}

\theoremstyle{remark}
\newtheorem{rem}[thm]{Remark}

\newcommand{\N}{\mathbb{N}}
\newcommand{\Z}{\mathbb{Z}}

\newcommand{\CA}{\mathcal{A}}
\newcommand{\CB}{\mathcal{B}}

\newcommand{\CE}{\mathcal{E}}

\newcommand{\CP}{\mathcal{P}}
\newcommand{\CQ}{\mathcal{Q}}

\newcommand{\CS}{\mathcal{S}}

\newcommand{\bs}\boldsymbol{}

\newcommand{\eq}[2]{ \begin{equation} \label{#1}\begin{split} #2 \end{split} \end{equation} }
\newcommand{\al}[1]{\begin{align} #1 \end{align} }
\newcommand{\als}[1]{\begin{align*} #1 \end{align*} }

\newcommand{\nn}{\nonumber \\}

\begin{document}

\title{When the sieve works}

\author{Andrew Granville}\thanks{Partially supported by NSERC}
\address{D\'epartement de math\'ematiques et de statistique\\
Universit\'e de Montr\'eal\\
CP 6128 succ. Centre-Ville\\
Montr\'eal, QC H3C 3J7\\
Canada}
\email{andrew@dms.umontreal.ca}

\author{Dimitris Koukoulopoulos}
\address{D\'epartement de math\'ematiques et de statistique\\
Universit\'e de Montr\'eal\\
CP 6128 succ. Centre-Ville\\
Montr\'eal, QC H3C 3J7\\
Canada}
\email{koukoulo@dms.umontreal.ca}

\author{Kaisa Matom\"aki}\thanks{Supported by the Academy of Finland grant no. 137883}
\address{Department of Mathematics\\
University of Turku\\
20014 Turku \\
Finland}
\email{ksmato@utu.fi}


\date{\today}

\begin{abstract} 
We are interested in classifying those sets of primes $\CP$ such that when we sieve out the integers up to $x$ by the primes in $\CP^c$ we are left with {\sl roughly} the expected number of unsieved integers. In particular, we obtain the first general results for sieving an interval of length $x$ with primes including some in $(\sqrt{x},x]$, using methods motivated by additive combinatorics.
\end{abstract}

\maketitle




\section{Introduction and motivation}

Let $\CE$ be a given subset of the primes $\leq x$.  The simplest sieve problem asks for estimates of $\CS(T, T+x; \CE)$, the number of integers $n$ in an interval $(T, T+x]$ which have no prime factors in the set $\CE$ (we write $(n,\CE)=1$ for convenience).
By a simple inclusion-exclusion argument one expects that the number of such integers is about
\[
x\prod_{p\in \CE}\left(1-\frac1p\right).
\]
This is provably always an upper bound, up to a constant:
\[
\CS(T, T+x; \CE) \ll x\prod_{p\in \CE}  \left( 1 -\frac 1p \right),
\]
and one gets the analogous lower bound
\[
\CS(T, T+x; \CE) \gg x\prod_{p\in \CE}  \left( 1 -\frac 1p \right),
\]
if $\CE$ is a subset of the primes up to $x^{1/2-o(1)}$ (see \cite[Theorem 11.13]{fi} noticing that the sieving limit $\beta = 2$ for $\kappa = 1$).  There seems to be little hope, in this generality, of increasing   the exponent ``$1/2$'' without major new ideas.  Moreover, one needs to have careful hypotheses: if, for instance,
$\CE=\{ p\leq \sqrt{x}\} \cup \{ x/\log x<p\leq x\}$ then
$\CS(0,x; \CE)=\pi(x/\log x)-\pi(\sqrt{x})+1\sim x/\log^2 x$, whereas our prediction was $\asymp x/\log x$.

In this article we will prove lower bounds on $\CS(0,x; \CE)$ in certain cases not covered by classical sieve theory, and use this evidence to guess at lower bounds in more generality. Let us first introduce some notation. Let $\CP$ be a given subset of the primes $\leq x$, and  $\CE=\{p\leq x:\ p\not\in \CP\}$, so that  $\CE \cup \CP$ is a partition of the primes $\leq x$. Let $\Psi(x;\CP)$ denote the number of integers up to $x$ all of whose prime factors are in $\CP$, that is
\[
\Psi(x;\CP):=\sum_{\substack{ n\leq x\\ p|n\implies p\in \CP}} 1  
 	= \sum_{\substack{ n\leq x\\ (n,\CE)=1}} 1 
	= \CS(0,x; \CE) .
\]
The inclusion-exclusion argument predicts that
\[
\Psi(x;\CP) \approx  x/u_{\CP} , 
\quad  \text{where} \quad
u_{\CP}:=\prod_{p\in \CE}  \left( 1 -\frac 1p \right)^{-1}.
\]
Hall \cite{hall} proved that $\Psi(x;\CP) \lesssim (e^\gamma/u_{\CP})\ x$, where $\gamma=0.5772156649\ldots$
is the Euler-Mascheroni constant. Subsequently, the authors of \cite{gs}  improved  this to
$\Psi(x;\CP) \lesssim ( e^\gamma/u_{\CP} - 1/{u_{\CP}^{2+o(1)}})   \ x$, and showed that this is ``best possible'' other than being more precise than the ``$o(1)$'' (which tends to $0$ when $u_{\CP} \to \infty$). However, in this paper, we are interested in lower bounds on $\Psi(x;\CP)$.

Hildebrand \cite{hild} showed that, $\Psi(x;\CP) \gtrsim \Psi(x; \CP')$ where $\CP'$ is the set of primes up to $x^{1/{u_{\CP}}}$ (note that here $u_{\CP'} \sim u_{\CP}$). Specifically, he showed that
\[
\frac{\Psi(x;\CP)}x\gtrsim\rho(u_{\CP});
\]
here $\rho(u)$ is the Dickman--de Bruijn function, defined by $\rho(u) =1$ for $0 \le u \le 1$, and $u\rho^{\prime}(u) = -\rho(u-1)$ for all $u \ge 1$.  One can show that $\rho(u) = 1/u^{u+o(u)}$, which is a lot smaller than the expected $1/u$. (See \cite{gs} for a different proof of Hildebrand's result.)

This last example is very special in that $\CP$ contains no large primes, and it is expected that other than for certain other extraordinary sets
$\CP$ one has
\[
\Psi(x;\CP)\asymp x\prod_{p\in \CE}  \left( 1 -\frac 1p \right).
\]
However this question has not really been studied in detail and there are other examples that must be taken into account:\ Let
\eq{comb-obstructions1}{
\CP=\bigcup_{1\leq m\leq N-1}\{ p \ \text{prime}:\ x^{m/(N+1)}<p<x^{m/N}\}.
}
Any product $\leq x$ of primes in these intervals lies in some interval of the same form. If now
$n=p_1p_2\dotsm p_k\le x$, where $x^{m_j/(N+1)}<p_j<x^{m_j/N}$, then $m_1+\dots +m_k\leq N$. The number of such integers $n$ with $m_1+\dots +m_k\leq N-1$ is $\leq x^{1-1/N}$. On the other hand, if $m_1+\dots +m_k=N$, then $k\geq 2$ (since each $m_j\leq N-1$) and therefore $\Psi(x;\CP) \ll_N x/(\log x)^2$, far smaller than the expected $x/u_{\CP}\asymp x/\log x$. A key thing to notice is that, in this example,
\[
\sum_{p\in \CP} \frac 1p \sim (N-1)\log(1+1/N)= 1 -\frac{3/2+o(1)}N  <1,
\]
by the prime number theorem, as $N\to\infty$. Hence we see that we can obtain a very small $\Psi(x;\CP)$ when $\sum_{p\in \CP} \frac 1p<1$.  One might guess that the issue in the last example is that there are no small primes in $\CP$. However, if we let
\eq{comb-obstructions2}{
\CP= \{ p \ \text{prime}:\ p\leq x^{1/N^2}\} \cup \bigcup_{1\leq m\leq N-1} \{ p \ \text{prime}:\  x^{m/(N+1)}<p<x^{m/N} \},
}
then we would also find that $\Psi(x;\CP)$ is far smaller than expected if $N\to\infty$ slowly enough.

One might also guess that the main issue in the above examples is that there are no large primes in $\CP$. However if we let
\[
\CP= \{ p \ \text{prime}:\ x^{1/v}<p\leq x^{1/u}\} \cup  \{ p \ \text{prime}:\   x^{1-1/v}<p\leq x \},
\]
then $\Psi(x;\CP) = \Psi(x;\CQ)+\pi(x)-\pi(x^{1-1/v})$ where $\CQ=\{ p:\ x^{1/v}<p\leq x^{1/u}\} $. Friedlander \cite{fried} established that $\Psi(x;\CQ)\sim \sigma(u,v) x/\log (x^{1/v})$ where $\sigma(u,v)=e^{-\gamma}\rho(u)+O(1/v\log v)$. Hence if $v>u^{(1+\epsilon)u}$  then
$\Psi(x;\CP)\sim e^{-\gamma}v\rho(u)  x/\log x$ as $u\to \infty$, whereas the ``expected'' value is
\[
x\prod_{p\in \CE}  \left( 1 -\frac 1p \right)\sim \frac{e^{-\gamma} v^2}{u(v-1)} \cdot \frac{x}{\log x} .
\]
Hence the ratio $\Psi(x;\CP)/ \text{``expected value''}$ is asymptotic to $u \rho(u)(1-1/v)$, which goes to $0$ rapidly  as $u\to \infty$.

So we see that the size of the primes in $\CP$ does not seem to  determine whether $\Psi(x;\CP)$ is close to its ``expected value''. Rather, we believe that the most important quantity in determining whether the sieve will work somewhat as expected is the largest $y$ for which  $\sum_{p\in \CP,\ p\geq y} \frac 1p>1$.


\begin{mainconj}\label{conj1} 
Fix $\epsilon>0$. There exists a positive constant $c$ such that if $\CP$ is a subset of the primes $\le x$ for which there is some $v\le c\sqrt{\log x}$ with
\[
\sum_{\substack{ p\in \CP \\ x^{1/ev}<p\leq x}} \frac 1p \geq 1+\epsilon,
\]
then
\[
\frac{\Psi(x;\CP)}x \ \ge\   A_{v}  \ \prod_{p\in \CE}  \left( 1 -\frac 1p \right),
\]
where $A_{v}$ is a constant with $A_{v}=v^{-v(1+o_\epsilon(1))}$ as $v\to\infty$.
\end{mainconj}




\begin{rem}\label{rem1.1}
Note that if $\CP=\{p\le x^{1/u}\}$, then $\Psi(x;\CP)/x\sim \rho(u)$. Also, taking $v \sim u$ in the above conjecture yields $\Psi(x;\CP)/x\gtrsim A_{v}/u$ and so $A_{v}\lesssim u\rho(u)=v^{-v(1+o(1))}$. One can even make the bolder guess that $A_{v}\sim v\rho(v)$.
\end{rem}


\begin{rem}\label{rem1.2}
Proposition \ref{prop2.3} below implies that if 
\eq{kappa}{
\kappa := \sum_{ \substack{  p\in \CP \\ x^{1/u}<p\leq x}} \frac{1}{p} 
}
is very small, then $\Psi(t;\CP)$ is indeed substantially smaller than expected for some $t\in [\sqrt{x},x]$. If, in addition, $\CP\subset[1,x^{1-\epsilon}]$, then we can take $t=x$. Hence it is certainly true that the size of $\kappa$ is important consideration. It would be interesting to obtain the strongest possible uniform version of Proposition \ref{prop2.3}.

On the other hand, Corollary \ref{cor2.2} below implies that as soon as $\kappa>\epsilon$ for some positive constant $\epsilon$, then there is a point $t \in [x^{1/u}, x]$ for which $\Psi(t;\CP)$ is of the expected size. It turns out to be a combinatorial problem to see if $t = x$ works and, in light of Bleichenbacher's theorem below, we believe that combinatorial obstructions, such the ones occurring for the sets $\CP$ given by \eqref{comb-obstructions1} and by \eqref{comb-obstructions2}, disappear as soon as $\kappa\ge1+\epsilon$.
\end{rem}


The main result in this paper is a weak form of Conjecture \ref{conj1}:


\begin{mainthm}\label{thm1}
There exist positive constants $\lambda$ and $c$ such that if $\CP$ is a subset of the primes $\le x$ for which there is some $v\le c\sqrt{\log x}$ with
\[
\sum_{\substack{ p\in \CP \\ x^{1/ev}<p\leq x}} \frac 1p \geq 1+\lambda,
\]
then
\[
\frac{\Psi(x;\CP)}x \ \gg \ \frac1{v^{O(v)}}  \ \prod_{p\in \CE}  \left( 1 -\frac 1p \right).
\]
\end{mainthm}


\begin{rem}\label{rem1.3}
One can superficially make Theorem \ref{thm1} appear stronger. For, if the ``$O(v)$'' in our lower bound is short for ``$Cv$'', then, for any $\epsilon \in (0, C)$, we can replace $\lambda$ in the hypothesis by $\Lambda:=\lambda+\ln(C/\epsilon)$ so that
\[
\sum_{\substack{ p\in \CP \\ x^{C/\epsilon ev}<p\leq x}} \frac 1p \geq
\sum_{\substack{ p\in \CP \\ x^{1/ev}<p\leq x}} \frac 1p  - \sum_{\substack{ p\in \CP \\x^{1/ev} <p\leq x^{C/\epsilon ev}}} \frac 1p \geq
1+\Lambda- \ln(C/\epsilon) = 1+\lambda,
\]
and hence Theorem \ref{thm1} implies, for $V=\epsilon v/C$, that
\[
\frac{\Psi(x;\CP)}x \ \gg \ \frac1{V^{CV}}  \ \prod_{p\in \CE}  \left( 1 -\frac 1p \right) \geq
\frac1{v^{\epsilon v}}  \ \prod_{p\in \CE}  \left( 1 -\frac 1p \right)  .
\]
\end{rem}


\begin{rem}\label{rem1.4}
It will be clear from the proof of Theorem \ref{thm1} that, under the same assumptions, one gets also the stronger conclusion
\[
\frac{\Psi(x;\CP) - \Psi(x-y; \CP)}{y} \ \gg \ \frac1{v^{O(v)}}  \ \prod_{p\in \CE}  \left( 1 -\frac 1p \right)
\]
for every $y > x^{1-1/(3ev)}$. Furthermore the conclusion of Theorem \ref{thm1} continues to hold if, more generally,
\[
\sum_{\substack{ p\in \CP \\ x^{1/ev}<p\leq x^{1/u}}} \frac 1p \geq \frac{1+\lambda}{u}
\]
for some $1 \leq u \leq v \leq c\sqrt{\log x}$, where
$\lambda$ and $c$ are absolute constants.
\end{rem}


The proof of our theorem involves a hodge-podge of techniques, from classical analytic number theory and sieve theory (as one might expect) to additive combinatorics, which seems to be new in this context (though \cite{lp} contains some related work as will be explained in Remark \ref{rem4.1}). Our starting point is the following result of Bleichenbacher \cite{bb} (see \cite[Section 9]{lp} for the proof) which may be viewed as a result in continuous additive combinatorics:


\theoremstyle{plain}
\newtheorem*{bbthm}{Bleichenbacher's Theorem}

\begin{bbthm}  If $u\geq 1$ and  $T$ is an open subset of $(0,\frac{1}{u})$ for which
\[
\int_{t\in T}\frac {dt}{t}\ > \ \frac 1u \ ,
\]
then there exist $t_1,t_2,\dots, t_k\in T$ for which $t_1+t_2+\dots+t_k=1$.
\end{bbthm}


Note that this is ``best possible'' since for the set $T_N=\bigcup_{j=1}^N (\frac j{N+1}, \frac jN)$ there is no solution to $t_1+t_2+\dots+t_k=1$ (as any sum of elements in intervals of the form $(\frac j{N+1}, \frac jN)$ is also in an interval of this form), yet
$\int_{T_N\cap(0,1/u)}\frac {dt}{t}\ \geq  [N/u] \log(1+1/N) > 1/u-2/N$, which tends to $1/u$ from below as $N\to \infty$. One  sees an analogy between this example and the first sieve example that we gave above.

The key idea in this paper is to reduce estimates for $\Psi(x;\CP)$ to quantitative questions of the type addressed in Bleichenbacher's theorem. Accordingly we make the following conjecture:


\theoremstyle{plain}
\newtheorem*{hypt}{Hypothesis T}

\begin{hypt} There exists a constant $\lambda_3>0$ such that if $1\le u\le v$ and $T$ is an open subset of $(\frac{1}{ev},\frac{1}{u})$ for which
\[
\int_{t\in T}\frac{dt}{t}\ge\frac {1+\lambda_3}{u} \ ,
\]
then there exists an integer $k\in[u,ev]$ and an absolute constant $\tau_v>0$  such that
\[
{\int \dots \int}_{\substack{ t_1+t_2+\dots+t_k=1\\t_1,t_2,\dots, t_k\in T}}\frac{dt_1dt_2\dotsm dt_{k-1}}{t_1t_2\cdots t_k}\ge \tau_v \left(  \int_{t\in T}\frac{dt}{t} \right)^k.
\]
\end{hypt}


By analogy with Bleichenbacher's Theorem, we conjecture that Hypothesis T holds for any choice of $\lambda_3>0$.  The importance of Hypothesis T can be seen in the following consequence, which will be proven in Section \ref{comb-equiv}.


\begin{mainprop}\label{prop1} If Hypothesis T is true for any fixed $\lambda_3>0$ and $\tau_v=v^{-O(v)}$, then
Conjecture \ref{conj1} holds for any fixed $\epsilon>0$ with $\kappa_v=v^{-O(v)}$.
\end{mainprop}


Actually we will formulate Hypotheses P and A which are analogous to Hypothesis T but Hypothesis P concerns counting primes $p_1, \dotsc, p_k \in \mathcal{P} \subseteq \mathbb{P} \cap [x^{1/ev}, x^{1/u}]$ for which $\log p_1 + \dotsc + \log p_k = x + O(1)$ and Hypothesis A concerns counting integers $a_1, \dotsc, a_k \in A \subseteq \mathbb{N} \cap (N/ev, N/u]$ for which $a_1 + \dotsc + a_k = N + O(k)$. We will show that all the three hypotheses are equivalent and that hypothesis P implies Conjecture~\ref{conj1}. 

We will use additive combinatorial tools to prove Hypothesis A for some sufficiently large counterpart of $\lambda_3$ (see Theorem \ref{thm7.1}) which allows us to deduce Theorem \ref{thm1} as desired. The value of $\lambda_3$ can be determined explicitly from the proof, but it will certainly not yield Hypothesis T  for {\it every} $\lambda_3>0$. See also Remark \ref{rem7.1} for more discussion about attainable $\lambda$-values.

Our results leave us wondering whether Conjecture \ref{conj1} and Theorem \ref{thm1} might be an indication of the truth for the more general problem of sieving intervals. Could it be that when we sieve an arbitrary interval of length $x$,  with a not-too-large subset of the primes up to $x$, then the number of integers left unsieved is predictable? There are only two types of examples known where we can do accurate calculations to better understand sieving: random intervals and intervals where most of the progressions sieved are $0 \pmod p$, and as far as we know the latter are where most extreme examples come from. Since we have now obtained some understanding of this usual source of extreme examples, we can speculate that this sort of criteria is the main issue, in general.

We conclude by mentioning that, aside from the theoretical interest of understanding the limitations of traditional heuristics in sieve methods, the generality of our results have applications beyond this subject. Indeed, in \cite{kaisa}, the third author discovered a rather unexpected application of these methods to counting real zeroes of holomorphic cusp forms.


\subsection*{Overview of the paper}\label{overview}

The paper is organized as follows: In Section \ref{log weights} we explore what happens when $\kappa$ defined by \eqref{kappa} is rather small. In particular, we show that as soon as $\kappa >\epsilon$, the quantity $\Psi(t;\CP)$ has the expected size for some $t\in[x^{1/u},x]$. Conversely, if $\kappa=o(1)$ as $u\to\infty$, then the size of $\Psi(t;\CP)$ is much smaller than expected for a certain $t\in[\sqrt{x},x]$. This makes it evident that in order for the sieve to work as expected, we need $\kappa$ to have some size. As Bleichenbacher's Theorem and the examples given in \eqref{comb-obstructions1} and in \eqref{comb-obstructions2} indicate, we should have that $\kappa>1$. However,  traditional sieve methods are incapable of handling this problem. Enter additive combinatorics. Indeed, as the results of Section \ref{technical} show, after some technical manipulations we can reduce the problem of bounding $\Psi(x;\CP)$ from below to counting $k$-tuples of primes $(p_1,\dots,p_k)\in (\CP\cap[x^{1/u},x])^k$ such that $\log p_1+\cdots +\log p_k=\log x+O(1)$, for some appropriate $k$. This reformulation of the problem, which we call {\it Hypothesis P} in analogy with Hypothesis T, makes clear the connection with additive combinatorics. In order to crystallize this connection even further and open the door to the use of additive combinatorial tools, in Section \ref{comb-equiv} we formulate the {\it Hypothesis A}, which is a purely combinatorial analogue of Hypothesis P and can be viewed as the discrete version of Hypothesis T. All these difference hypotheses are, in fact, equivalent as we show in Section \ref{comb-equiv}. It is Hypothesis A that we will eventually prove in Section \ref{hypa-proof}, using some tools of discrete additive combinatorics developed in Section \ref{add-comb}. Finally, in Section \ref{a-p-t}, we explore further the connections between our three difference hypotheses, A, P and T. 


\section{Sieving with logarithmic weights}\label{log weights}

If we introduce the weight $1/n$ at each integer $n$ (the so-called ``logarithmic weights''), then we simplify the problem enormously:


\begin{lem}\label{lem2.1} 
If $\CP$ is a subset of the primes $\le x$ and $\CE=\{p\le x\}\setminus\CP$, then
\[
\prod_{p\in \CE} \left( 1-\frac 1p\right)\lesssim \frac 1{\log x} \ \sum_{\substack{ n\leq x \\ p|n \implies p\in \CP}} \frac 1n \lesssim
e^\gamma \prod_{p\in \CE} \left( 1-\frac 1p\right),
\]
where $\gamma$ is the Euler-Mascheroni constant.
\end{lem}


\begin{proof}
Let us first prove the lower bound. We have
\[
\prod_{p\in \CE} \left( 1-\frac 1p\right)^{-1} \geq \sum_{\substack{ \ell\leq x \\ p|\ell \implies p\in \CE}} \frac 1\ell,
\]
so that
\[
\sum_{\substack{ m\leq x \\ p|m \implies p\in \CP}} \frac 1m
\geq \prod_{p\in \CE} \left( 1-\frac 1p\right)\sum_{\substack{ \ell\leq x \\ p|\ell \implies p\in \CE}} \frac 1\ell
\sum_{\substack{ m\leq x \\ p|m \implies p\in \CP}} \frac 1m
\geq \prod_{p\in \CE} \left( 1-\frac 1p\right)\sum_{\substack{ n\leq x  }} \frac 1n,
\]
since every integer $n\leq x$ may be written as $\ell m$.

On the other hand we have the upper bound
\[
\sum_{\substack{ n\le x\\p|n\implies p\in\CP}} \frac1n\le\prod_{p\in \CP} \left( 1-\frac 1p\right)^{-1} = \prod_{p\leq x} \left( 1-\frac 1p\right)^{-1}\prod_{p \in \CE} \left( 1-\frac 1p\right)\sim e^\gamma\log x\prod_{p\in\CE}\left(1-\frac1p\right),
\]
by Mertens' theorem.
\end{proof}


\begin{rem}
Note that
\eq{eq2.1}{
\sum_{\substack{ n\leq x \\ p|n \implies p\in \CP}} \frac 1n = \frac{\Psi(x;\CP)}x + \int_1^x \frac{\Psi(t;\CP)}{t^2} dt,
}
so that Lemma \ref{lem2.1} can be re-phrased as a weighted mean of $\Psi(t;\CP)$-values:
\[
\prod_{p\in \CE} \left( 1-\frac 1p\right) \lesssim
\int_1^x \frac{\Psi(t;\CP)}{t} \cdot \frac{dt}{t} \bigg/ \int_1^x   \frac{dt}{t}
\lesssim e^\gamma \prod_{p\in \CE} \left( 1-\frac 1p\right),
\]
since $\Psi(x;\CP)/x \ll \prod_{p\in \CE} \left(1-\frac1p\right) $ by classical sieve theory (as discussed in the introduction).
\end{rem}


We can use Lemma \ref{lem2.1} to prove a first lower bound in the direction of Theorem \ref{thm1}, though with a different emphasis: We show that as soon as $\sum_{x^{1/u}<p\le x, p \in \CP}\frac 1p$ is at least $\epsilon$, there is some $t \in [x^{1/u}, x]$ for which $\Psi(t;\CP)$ is of expected size.


\begin{cor}\label{cor2.2}  
Fix $\epsilon\in(0,1)$. Suppose that $\CP$ is a subset of the primes $\le x$, and $u\in[1,\log x]$ is such that
\[
\sum_{\substack{ p\in \CP \\ x^{1/u}<p\leq x}} \frac 1p>\epsilon.
\]
Then there exists $t\in [x^{1/u},x]$ such that
\eq{eq2.2}{
\frac{\Psi(t;\CP)}t \gg \frac {\epsilon\min\{1,\epsilon u\}}{\log u} \ \prod_{\substack{ p\in \CE \\ p\leq t}}
		 \left( 1 -\frac 1p \right).
}
\end{cor}


\begin{proof}
By \eqref{eq2.1}
\[
\int_{x^{1/u}}^x\frac{\Psi(t;\CP)}{t^2}dt \geq \sum _{\substack{ x^{1/u} \leq n \leq x \\ p|n \implies p \in \CP }} \left(\frac{1}{n}-\frac{1}{x}\right) \geq \frac{1}{2} \sum _{\substack{ x^{1/u} \leq n \leq x/2 \\ p|n \implies p \in \CP }} \frac{1}{n}.
\]
Writing here $n = ab$, where prime factors of $a$ are $\leq x^{1/u}$ and prime factors of $b$ are $> x^{1/u}$, and discarding some $n$, we see that
\[
\int_{x^{1/u}}^x\frac{\Psi(t;\CP)}{t^2}dt
\ge\frac12 \sum_{\substack{ a\le x^{\epsilon/2}/2\\ p|a \implies p \in \CP \cap [1, x^{1/u}]}} \frac1a \sum _{\substack{ 1<b\le x^{1-\epsilon/2} \\ p|b \implies p \in \CP \cap (x^{1/u}, x^{1-\epsilon/2}] }} \frac 1b.
\]
Here
\[
S:=\sum _{\substack{ 1<b\le x^{1-\epsilon/2} \\ p|b\implies p\in \CP \cap (x^{1/u}, x^{1-\epsilon/2}] }} \frac 1b\ge\sum _{\substack{ p\in\CP \\ x^{1/u}<p\le x^{1-\epsilon/2}}} \frac1p
>\epsilon+\log(1-\epsilon/2)+o(1)\ge\frac\epsilon4,
\]
hence $S\geq \frac \epsilon 8(1+S)$ and so
\[
\int_{x^{1/u}}^x\frac{\Psi(t;\CP)}{t^2}dt \geq \frac \epsilon 8 \ \sum_{\substack{ a\le x^{\epsilon/2}/2\\p|a \implies p \in \CP \cap [1, x^{1/u}] }} \frac1a \sum _{\substack{ b\le x^{1-\epsilon/2} \\ p|b \implies p \in \CP \cap (x^{1/u}, x^{1-\epsilon/2}] }} \frac 1b.
\]
Consequently, Lemma \ref{lem2.1} implies that
\als{
\int_{x^{1/u}}^x\frac{\Psi(t;\CP)}{t^2}dt 
	&\gg \epsilon \prod_{\substack{ p\le x^{\epsilon/2}\\ p\in \CP \cap [1,x^{1/u}]}}
		 \left(1+ \frac1p \right) \prod_{p\in\CP\cap(x^{1/u},x^{1-\epsilon/2}]} \left(1+\frac1p\right)\\
	&\gg \epsilon\min\{1,\epsilon u\} \prod_{p\in\CP} \left(1+\frac1p\right)
		\asymp \epsilon\min\{1,\epsilon u\}\log x \prod_{p\in\CE} \left(1-\frac1p\right).
}
If now for every $t\in[x^{1/u},x]$ we have that
\als{
\frac{\Psi(t;\CP)}t   
	\le \eta\cdot \frac{\epsilon \min\{1,\epsilon u\}}{\log u}
		 \prod _{\substack{ p\in\CE\\p\le t}} \left(1-\frac1p \right)
	\ll \eta\cdot\frac {\epsilon \min\{1,\epsilon u\}}{\log u} \frac{\log x}{\log t} 
		\prod_{p\in\CE} \left(1-\frac1p \right),
}
then
\als{
\int_{x^{1/u}}^x\frac{\Psi(t;\CP)}{t^2} 
	&\ll \eta\cdot\frac {\epsilon \min\{1,\epsilon u\}}{\log u} \log x 
		\prod_{p\in\CE} \left(1-\frac1p \right) \int_{x^{1/u}}^x \frac{dt}{t\log t} \\
	&\ll \eta\cdot \epsilon \min\{1,\epsilon u\}\log x\prod_{p\in\CE} \left(1-\frac1p \right).
}
Choosing $\eta$ small enough, we arrive at a contradiction. So the claimed result follows.
\end{proof}


\begin{rem}\label{rem2.2} 
Let $T$ be the set of  $t\in [x^{1/u},x]$ for which \eqref{eq2.2} holds. By the same proof, and the usual sieve upper bound, we obtain
$\int_{t\in T} dt/(t\log t) \gg \epsilon\min\{1,\epsilon u\}$
\end{rem}


The lower bound \eqref{eq2.2} obtained here is much better than the lower bound in Theorem \ref{thm1}, but it only works for some values of $t$. One cannot essentially improve Corollary \ref{cor2.2} in general, at least when $1/\epsilon \ll u \ll \sqrt{\log x}$: Take a prime $q \asymp (\log u) / \epsilon$ and let $\CP$ be the set of primes which are $\equiv 1 \pmod{q}$. Then the classical sieve yields
\als{
\frac{\Psi(t, \CP)}{t} 
	&= \frac{1}{t}|\{n = qk+1 \leq t \colon p \mid n \implies p \not \in \CE \setminus \{q\} \}| \\
	&\ll \frac{1}{q} \prod_{_{\substack{ p \in \CE \setminus\{q\} \\ p \leq (t/q)^{1/2} }}} 
		\left(1-\frac{1}{p}\right) \asymp \frac{\epsilon}{\log u} 
		\prod_{_{\substack{ p \in \CE \\ p \leq t }}} 
		\left(1-\frac{1}{p}\right)
}
for every $t \in [x^{1/u}, x]$.

Next, we prove a converse result to Corollary \ref{cor2.2}, but first we need an estimate which belongs to the theory of smooth numbers. Its proof is an application of Rankin's method, together with an additional averaging which recovers a logarithmic loss that occurs in the original version of Rankin's method. It can be found, for example, in Kevin Ford's notes \cite[Theorem $\Psi$]{kf}, though it is possible that it has appeared before in the literature. We give the full proof for completeness.


\begin{prop}\label{smooth}
Let $x\ge3$ and $u\ge1$ such that $u\le (1/2-\epsilon)\log x/\log\log x$, for some fixed $\epsilon\in(0,1/3)$. If $\CP$ is a subset of the primes $\le x^{1/u}$, then
\[
\frac{\Psi(x;\CP)}{x} \ll_\epsilon \frac{e^{O(u)}}{(u\log u)^u} \prod_{\substack{ p\in \CE \\ p\leq x}} \left( 1 -\frac 1p \right).
\]
\end{prop}


\begin{proof} Without loss of generality, we may assume that $u$ is large enough. Set $y=x^{1/u}$ and note that $u\le y^{1/2-\epsilon}/\log y$ by our assumption that $u\le (1/2-\epsilon)\log x/\log\log x$. In particular, we may assume that $y$ is large too. 

Our starting point is the identity
\eq{rankin-1}{
\sum_{ \substack{ n\le x \\ p|n\implies p\in\CP}} \log n 
	= \sum_{\substack{m\le x \\ p|m\implies p\in\CP }}\ 
		\sum_{\substack{d \le x/m \\ p|d\implies p\in\CP}}\Lambda(d) .
}
Fix some $\delta\in [1/\log y,1/2-\epsilon]$ and note that, for $1\le n\le x$,
\[
\log x = \log n + \log(x/n) \le \log n + \frac1{1-\delta} \cdot \frac{x^{1-\delta}}{n^{1-\delta}}
	\le \log n + \frac{6x^{1-\delta}}{n^{1-\delta}} .
\]
Together with \eqref{rankin-1}, this implies that
\eq{rankin-2}{
(\log x)\Psi(x;\CP)  
	\ll \sum_{ \substack{ n\le x \\ p|n\implies p\in\CP }}  \frac{x^{1-\delta}}{n^{1-\delta}}
 		+  \sum_{\substack{m\le x \\ p|m \implies p\in\CP\ }}\ 
		\sum_{\substack{d \le x/m \\ p|d\implies p\in\CP }}\Lambda(d) .
}
Next, note that 
\[
\sum_{\substack{d \le x/m \\ p|d\implies p\in\CP  }}\Lambda(d)
	\le \sum_{p\le \min\{y,x/m\}} (\log p) \sum_{\substack{ \nu\ge1 \\ p^\nu\le x/m }}1
	\ll \sum_{p\le \min\{y,x/m\} } \log(x/m),
\]
by our assumption that $\CP\subset\{p\le y\}$. So, if $x/y<m\le x$, then we find that
\eq{rankin-3}{
\sum_{\substack{d \le x/m \\ p|d \implies p\in\CP }}\Lambda(d)
	\ll\frac{x}{m}
	\le \frac{y^\delta x^{1-\delta}}{m^{1-\delta}},
}
whereas, if $1\le m\le x/y$, then
\als{
\sum_{\substack{d \le x/m \\ p|d \implies p\in\CP}}\Lambda(d)
	\ll \frac{y}{\log y} \log (x/m) \ll \frac{y^\delta (x/m)^{1-\delta}}{\log (y^\delta (x/m)^{1-\delta})} \log (x/m) \ll \frac{y^\delta x^{1-\delta}}{m^{1-\delta}}.
}
In any case, the estimate \eqref{rankin-3} does hold. 
%
Combining it with\ \eqref{rankin-2}, we deduce that
\als{
\Psi(x;\CP)
	\ll \frac{y^\delta x^{1-\delta}}{\log x} \sum_{ \substack{ n\le x \\ p|n\implies p\in\CP}}\frac{1}{n^{1-\delta}}
	&\le \frac{y^\delta x^{1 - \delta}}{\log x} \prod_{p\in\CP } \left(1 - \frac1{p^{1 - \delta}} \right)^{-1} \\
	&\ll_\epsilon y^\delta x^{1 - \delta}\exp \left\{ \sum_{p\in\CP} \frac{p^\delta - 1}{p} \right\}  	
		 \prod_{\substack{ p\in\CE   \\  p\le x}} \left(1-\frac{1}{p}\right)  ,
}
by our assumption that $\delta\le 1/2-\epsilon$. Finally, note that
\als{
\sum_{p\in \CP }    \frac{p^{\delta} - 1}{p} 
	\le \sum_{p\leq y }    \frac{p^{\delta} - 1}{p} 	
	\ll   \frac{e^{\delta\log y}}{\delta\log y}
}
by the Brun-Titchmarsch inequality. So writing $\delta=v/\log y$, we arrive to the estimate
\[
\frac{\Psi(x;\CP)}{x} \ll    \frac{e^{O(e^v/v)} } {e^{uv}} .
\]
We choose $v$ such that $e^{v}/v=u$. This produces a value of $\delta$ in the interval $[1/\log y,1/2-\epsilon]$ as long as $u=e^{\delta\log y}/(\delta\log y)$ is in the interval $[e,y^{1/2-\epsilon}/((1/2-\epsilon)\log y)]$, which does hold. Since $v=\log u + \log\log u+O(1)$, the proposition follows.
\end{proof}


\begin{prop}\label{prop2.3}
Suppose that $\CP$ is a subset of the primes $\le x$, and $u\in[1,\log x]$ is such that
\[
\sum_{\substack{ p\in \CP \\ x^{1/u}<p\leq x}}  \frac 1p :=\kappa\ll 1.
\]
There exists $t\in [x^{1/2},x]$ such that
\[
\frac{\Psi(t;\CP)}t \ll  \left( \kappa +
x^{-1/6}+u^{-u/4} \right) \ \prod_{\substack{ p\in \CE \\ p\leq t}} \left( 1 -\frac 1p \right).
\]
If, in addition, $\CP\subset[1,x^{1-\epsilon}]$, then we can take $t=x$ provided we replace $\kappa$ by $\kappa/\epsilon$. In either case, if $\kappa=o(1)$ as $u\to\infty$ and $x\to\infty$, then $\Psi(t;\CP)/t$ is much smaller than expected.
\end{prop}


\begin{proof} First, we show the second claim because its proof is simpler. Note that
\[
\Psi(x;\CP) \le \Psi(x;\CP\cap[1,x^{1/u}]) 
	+ \sum_{\substack{ n\le x \\ p|n \implies p\in\CP}} \sum_{\substack{p|n \\ x^{1/u}<p\le x}} 1 .
\]
The first sum is $\ll x(u^{-u} + x^{-1/3}) \prod_{p\in\CE\cap[1,x]}(1-1/p)$, by Proposition \ref{smooth} applied with $\min\{u,(\log x)/(2.5\log\log x)\}$ in place of $u$. The second sum equals
\[
\sum_{\substack{p\in\CP \\ x^{1/u}<p\le x}} \sum_{\substack{ m\le x/p \\ p|m\implies p\in\CP}} 1
	\ll \sum_{\substack{p\in\CP \\ x^{1/u}<p\le x}} \frac{x}{p\log(x/p)} \prod_{p\in\CP\cap[1,x/p]} \left(1+ \frac{1}{p} \right) 
	\ll \frac{\kappa x}{\epsilon\log x} \prod_{p\in\CP} \left(1+ \frac{1}{p} \right) ,
\]
by our assumption that $\CP\subset[1,x^{1-\epsilon}]$. Therefore
\[
\frac{\Psi(x;\CP)}{x} \ll  \left( u^{-u} + x^{-1/3} + \kappa/ \epsilon \right) \prod_{\substack{p\in\CE \\ p\le x}} \left(1 - \frac{1}{p} \right),
\]
as claimed.

Finally, we show the first part of the proposition. We may assume that $u$ is large enough and $u = o(\log x)$. Our starting point is the relation
\eq{kappa small e1}{
 \int_{\sqrt{x}}^x \frac{\Psi(t;\CP)}{t} \frac{dt}{t} + \frac{\Psi(x;\CP)}{x} 
 	=   \sum_{\substack{ \sqrt{x}<n\leq x \\ p|n\implies p\in \CP  }} \frac 1n +  \frac{\Psi(\sqrt{x};\CP)}{\sqrt{x}} ,
}
which follows by integration by parts. If we show that each term on the right hand side of \eqref{kappa small e1} is 
\[
\ll \prod_{\substack{ p\le x  \\ p\in \CE}} \left( 1- \frac{1}{p} \right) (\kappa + x^{-1/6} + u^{-u/4} ) \log x,
\]
then the claimed result follows, by taking the minimum of $\Psi(t;\CP)/t$ on the left side of \eqref{kappa small e1}.

First, we bound the sum over $n$. Let $v=\min\{u,(\log x)/(2.5\log\log x)\}$ and set $y=x^{1/v}\ge \max\{x^{1/u},(\log x)^{2.5}\}$. Notice that
\[
\sum _{\substack{ b \geq 1 \\ p \mid b\implies p \in \CP \cap (y, x] }} \frac 1b \leq \prod_{\substack{ y<p\leq x \\ p\in \CP }} \left( 1 - \frac 1p \right)^{-1} \leq e^{O(\kappa)} \ll 1.
\]
So, writing $n = ab$ with $a$ having prime factors $\leq y$ and $b$ having prime factors $>y$, and adding some extra $ab$, we see that
\al{
\sum_{\substack{ \sqrt{x}<n\leq x \\ p|n\implies p\in \CP  }} \frac 1n 
	&\leq \sum_{\substack{ a\leq \sqrt{x} \\ p|a\implies p \in \CP \cap [1, y] }} \frac 1a
		\cdot \sum_{\substack{ b>1 \\ p|b\implies p \in \CP \cap (y, x] }} \frac 1b +
		\sum_{\substack{ \sqrt{x}<a\leq x \\ p|a\implies p \in \CP \cap [1, y] }} \frac 1a
		\sum_{\substack{ b\geq 1 \\ p|b\implies p \in \CP \cap (y, x] }} \frac 1b \nn
	&\leq \sum_{\substack{ a\geq 1 \\ p|a\implies p \in \CP \cap [1, y] }} \frac 1a
		\cdot  ( e^{O(\kappa)}-1) +
		\sum_{\substack{ \sqrt{x}<a\leq x \\ p|a\implies p \in \CP \cap [1, y] }} \frac 1a
		\cdot  e^{O(\kappa)} \nn
	 &\ll \kappa \prod_{\substack{ p\leq y \\ p\in \CP }} \left( 1 - \frac 1p \right)^{-1} 
		+ (\log x) \frac{e^{O(v)}}{(v\log v)^{v/2}}  \prod_{\substack{ p\leq x \\ p\in \CE }} \left( 1 - \frac 1p \right)  \nn
	&\ll (\log x) (\kappa + x^{-1/3} + (2u)^{-u/2} ) \prod_{\substack{ p\leq x \\ p\in \CE }} \left( 1 - \frac 1p \right) ,   \label{kappa small e2}
}
by Proposition \ref{smooth} and partial summation. We deduce that
\[
\Psi(x;\mathcal{P})\leq \sqrt{x}+\sum_{\substack{\sqrt{x}<n\leq x \\ p|n\implies p\in \mathcal{P}}} \frac xn \ll  x \prod_{\substack{ p\leq x \\ p\in \mathcal{E}}} \left( 1 - \frac 1p \right)  \left( \kappa + x^{-1/3}+(2u)^{-u/2}\right) \log x.
\]
Since $(\sqrt{x}^{1/(u/2)},\sqrt{x}]\subset(x^{1/u},x]$, applying the above relation with $x$ and $u$ replaced by $\sqrt{x}$ and $u/2$, respectively, we see that
\eq{kappa small e3}{
\frac{\Psi(\sqrt{x};\CP)}{\sqrt{x}} \ll  \prod_{\substack{ p\leq x \\ p\in \CE }} \left( 1 - \frac 1p \right)  \left( \kappa + x^{-1/6}+u^{-u/4}\right) \log x.
}
Inserting relations \eqref{kappa small e2} and \eqref{kappa small e3} into \eqref{kappa small e1} completes the proof of the proposition.
\end{proof}


This last estimate is much smaller than one might have guessed given Lemma \ref{lem2.1}.

We have now seen that if there are very few large primes in $\CP$ then one can improve the sieve upper bounds for some values of $t$. Finally, note that the assumption that $\CP\subset[1,x^{1-\epsilon}]$ is essential in the second part of Proposition \ref{prop2.3}. Indeed, if $\CP=\{ p\leq x^\epsilon\} \cup \{ p:\ x^{1-\epsilon}<p\leq x\}$, then $\Psi(x,\CP)$ contains all integers of the form $mp$ where $p$ is a prime in the range $x^{1-\epsilon}<p\leq x$, and $m\leq x/p$. Hence
\[
\Psi(x;\CP) \geq \sum_{x^{1-\epsilon}<p\leq x} \frac xp \gg \epsilon x \asymp x \prod_{\substack{ p\in \CE \\ p\leq x}} \left( 1 -\frac 1p \right) .
\]


\section{Technical reductions}\label{technical}

The hypothesis in Theorem \ref{thm1} relies on there being a reasonable density of  ``large'' primes in $\CP$. More generally, we may ask what happens when some interval $(x^{1/ev},x^{1/u}]$ contains lots of primes of $\CP$. Reducing to the analogous problem where $\CP$ is now restricted to be a subset of the primes in this interval, we formulate


\theoremstyle{plain}
\newtheorem*{hypp}{Hypothesis P}

\begin{hypp}
There exist constants $\lambda_1>0$ and $C_1 > 1$ such that if $v^2\le \lambda_1 \log x/C_1$, $1\le u\le v$ and $\CP$ is a subset of the primes in $(x^{1/ev},x^{1/u}]$ for which
\eq{eq4.1}{
\sum_{\substack{ p\in \CP}}\frac 1p \geq \frac {1+\lambda_1}u ,
}
then for any $\delta \in [x^{-1/(3ev)}, 1/2]$ there exists an integer $k\in[u,ev]$ and an absolute constant $\pi_v>0$ such that
\[
\sum_{\substack{  (p_1,\dotsc,p_k)\in \CP^k\\  (1-\delta)x < p_1\cdots p_k<x}} \frac 1{p_1\dotsm p_k} \ge \frac{\delta \pi_v}{\log x} \cdot \left( \sum_{\substack{ p\in\CP}}\frac1p\right)^k\ .
\]
\end{hypp}


\begin{rem}\label{rem4.1}
In \cite{lp} Bleichenbacher's theorem is used to prove a result like Hypothesis P but with a logarithmic loss $(\log x)^{-O(v)}$ in the obtained lower bound --- see Theorem 4 and Proposition 10.1 there.
\end{rem}


\begin{proof}[Proof that Hypothesis P with $\pi_v=v^{-O(v)}$ implies Theorem \ref{thm1} with $\lambda=\lambda_1+\epsilon$ and $c = \sqrt{\lambda_1/(2C_1)}$] 
Set $\eta=\min\{\epsilon,1\}/3$. Let $\CA=\CP\cap[1,x^{1/ev}]$ and $\CB=\CP\cap(x^{1/ev},x]$ so that
\[
\Psi(x;\CP)\geq\sum_{\substack{ a\leq x^\eta\\p|a\implies p\in \CA}}\Psi(x/a;\CB),
\]
since we can write any $n$ composed only of prime factors from $\CP$ as $n=ab$ where $a$ and $b$ are composed only of prime factors from $\CA$ and $\CB$, respectively.
For each  $a\le x^\eta$, we have that
\[
\sum_{\substack{ p\in \CB \\(x/a)^{1/ev}<p\leq x/a}}\frac 1p\ge\sum_{\substack{ p\in \CP \\x^{1/ev}<p\leq x}}\frac 1p+\log(1-\eta)+o(1)\ge1+\lambda_1+\epsilon-\frac{\eta}{1-\eta}+o(1)>1+\lambda_1,
\]
as $x\to\infty$. Applying Hypothesis P with $u=1$ and  $\delta=1/2$ to the set $\CB$ yields
\[
\Psi(x/a;\CB) \ge \sum_{\substack{ x/2a<n\leq x/a \\ p|n\implies p\in \CB}} \frac{x/2a}{n}
\gg\frac1{v^{O(v)}}\,\cdot\,\frac x{a\log x} \,
\]
and consequently
\[
\frac{\Psi(x;\CP)}x\gg\frac1{v^{O(v)}}\frac1{\log x}\sum_{\substack{ a\leq x^\eta\\p|a\implies p\in \CA}}\frac1a
\gg\frac{\eta}{v^{O(v)}}\prod _{\substack{ p\le x^\eta \\ p\in\CE }} \left(1-\frac1p\right)
\gg\frac{\eta}{v^{O(v)}}\prod_{p\in\CE} \left(1-\frac1p\right),
\]
by Lemma \ref{lem2.1}, which completes the proof of Theorem \ref{thm1}.
\end{proof}


So we have shown that in order to prove Theorem \ref{thm1} it suffices to prove the more convenient Hypothesis P with certain choices of the parameters therein.


\section{Equivalent problems in combinatorics}\label{comb-equiv}

In our sieve question we are seeking to sieve the integers up to $x$ by a given set of primes $\CE$ which is, as discussed in the introduction, the same thing as  counting the number of integers up to $x$ that are composed of primes from a given set $\CP$. This makes this a rather special case of sieving an interval, since the problem can now be approached as a question of counting lattice points:\
If $p_1 \cdots p_k \leq x$, then
\[
 \log p_1 + \log p_2 +\cdots + \log p_k\leq \log x
\]
and there are various techniques for attacking this problem. However they are not really effective, since here we have an enormous dimension compared to the volume of our region, even when restricting the primes in $\CP$ to an interval $[x^{1/ev}, x^{1/u})$. We can however cut the dimension of the problem significantly by taking approximations that do not greatly effect the answer. For example, if we replace each $\log p$ by $[\log p]$ and take $N$ to be an integer close to $\log x$, then we can count integer solutions to $a_1+a_2+\cdots + a_m\leq N$, and weight each $a_i$ by the number of primes $p$ in $\CP$ for which $[\log p]=a$.  However even this problem is of rather high dimension to directly use lattice point counting results, so instead we attack this as a question in combinatorics.


\theoremstyle{plain}
\newtheorem*{hypa}{Hypothesis A}

\begin{hypa}
There exist constants $\lambda_2>0$ and $C_2>1$ such that if $v^2\le \lambda_2N/C_2$, $1\le u\le v$ and $A$ is a subset of the integers in $(\frac{N}{ev},\frac{N}{u}]$ such that
\[
\sum_{\substack{ a\in A  }}  \frac 1a \geq \frac {1+\lambda_2}{u} ,
\]
then there exists an integer $k\in[u,ev]$, an absolute constant $\alpha_v>0$ and an integer $n\in [N-k,N]$ such that
\[
\sum_{\substack{  (a_1,\dots,a_k)\in A^k\\  a_1+\dots+a_k=n}} \frac 1{a_1\dotsm a_k}
\ge\frac{\alpha_v}{N} \left(  \sum_{\substack{ a\in A  }}  \frac 1a \right)^k .
\]
\end{hypa}


\begin{prop}\label{prop5.1}
{\rm (i)} Hypotheses P and A are equivalent, with $\lambda_2\asymp \lambda_1$ and  
$\alpha_v e^{O(v)} \gg \pi_v\gg \alpha_v\min\{1,\lambda_1^v\}v^{-O(v)}$.
\end{prop}


We will prove this at the end of this section.  We first note reasons for some of the conditions in Hypothesis A:

\begin{itemize}
\item
If $A$ is the set of integers in $(N/(k+1), N/k-1)$ then there are no sums of elements of $A$ in the interval $[N-k, N]$ and $\sum_{n\in A}1/n\sim\log(1+1/k)=1/k+O(1/k^2)$. Hence we must have $\lambda_2\geq 0$. However, we do believe that Hypothesis A holds for any $\lambda_2>0$.

\item
If $A$ is the set of integers $\equiv 0 \pmod d$ in $(\frac{N}{ev},\frac{N}{v}]$, then
$\sum_{a\in A}  1/a \sim 1/d$ as $N\to\infty$, and there are no solutions to $a_1+\cdots +a_k=n$ for any $n$ in an interval $md<n<(m+1)d$. Hence $n$ must be chosen from an interval of length $\geq d$. This explains the length of the interval for $n$ in Hypothesis A.


\end{itemize}

We will eventually prove Hypothesis A with $\alpha_v=1/v^{O(v)}$ and big enough constant $\lambda_2$. By Proposition \ref{prop5.1}(i), this yields Hypothesis P with $\pi_v=1/v^{O(v)}$, and therefore Theorem \ref{thm1} as was shown in Section \ref{technical}.

Since Hypothesis A involves so many integers, one might think to approximate the set of integers $A$ by a continuous variable; for instance, by considering very short intervals around each $a/N$, so as to obtain Hypothesis T (which is stated in the introduction):


\theoremstyle{plain}
\newtheorem*{prop 5.1(ii)}{Proposition \ref{prop5.1}}

\begin{prop 5.1(ii)}
{\rm (ii)} Hypotheses A and T are equivalent, with $\lambda_2\asymp \lambda_3$ and $\tau_ve^{O(v)}\gg \alpha_v\gg \tau_v/e^{O(v)}$.
\end{prop 5.1(ii)}

Combining this with Proposition \ref{prop5.1}(i) and the result from Section \ref{technical}, we can deduce Proposition \ref{prop1}. 







\subsection*{Proof of Proposition \ref{prop5.1}}
We conclude this section with the proof that our three hypotheses, A, P and T, are equivalent.


\begin{proof}[Proof that Hypothesis A  implies Hypothesis T] 
Assume that $T$ is an open subset of $(1/ev,1/u)$ such that $\int_Tdt/t \ge(1+\lambda_3)/u$. An open subset of the reals is a union of disjoint open intervals, and the number of intervals in the union is countable (as may be seen by labelling each interval with some rational it contains). Hence we may write $T=\bigcup_{i\ge1}(\alpha_i,\beta_i)$. But then there exists an integer $m$ such that if $S=\bigcup_{i=1}^m (\alpha_i,\beta_i)$, then $\int_{t\in S} dt/t\ge (1+2\lambda_3/3)/u$. By replacing $T$ in our assumption by $S$ and $\lambda_3$ by $2\lambda_3/3$, we may assume that $T$ is a finite union of open intervals.  We select $N$ to be much larger than $\max_i \{ v/|\alpha_i-\beta_i|\}$, $mv^3/\lambda_3$ and $mv^4$. Let
\[
A=\bigcup_{i=1}^m  \{a\in \Z:\ \alpha_iN+2ev <a < \beta_iN-2ev\}.
\]
Since $\int_{t=a/N}^{(a+1)/N} dt / t = 1/a+O(1/a^2)$, we deduce that
\[\sum_{a\in A}\frac1a= \int_{t\in T}\frac{dt}t+O\left(\frac {mv^2}{N}\right)=\left(1+O\left(\frac{mv^3}N\right)\right)\int_{t\in T}\frac{dt}t \ge\frac{1+\lambda_3/2}u.\]
Now, if $\lambda_3\ge2\lambda_2$, then Hypothesis A implies that there exists an integer $k\in[u,ev]$ and an integer $n\in[N-k,N]$ such that
\[
\sum_{\substack{  (a_1,\dots,a_k)\in A^k\\  a_1+\dots+a_k=n}} \frac 1{a_1\dotsm a_k}
\ge \frac{\alpha_v}N \left(  \sum_{\substack{ a\in A  }}  \frac 1a \right)^k \geq \frac{\alpha_v e^{-O(k)}}N \left(  \int_{t \in T} \frac {dt}{t} \right)^k.
\]
For each $k$-tuple $(a_1,\dots,a_k)\in A^k$ with $a_1+\cdots+a_k=n$, consider $t_i\in (a_i/N, (a_i+1)/N)\subset T$ for $1\leq i\leq k-1$ and define $t_k=1-(t_1+\cdots+t_{k-1})$. Then we have that
\[
\left|t_k-\frac{a_k}N\right| = \left|\frac{N-n}{N} + \sum_{i=1}^{k-1} \left(\frac{a_i}{N}-t_i\right)\right| \le\frac{2k-1}{N}\le\frac{2ev}{N}
\]
and consequently $t_k\in T$. Hence
\als{
{\int \dots \int}_{\substack{   t_1+t_2+\dots+t_k=1\\t_1,t_2,\dots, t_k\in T}} 
		\frac{dt_1dt_2\dotsm dt_{k-1}}{t_1t_2\cdots t_k}
			&\gg \sum _{\substack{  (a_1,\dots,a_k)\in A^k\\  a_1+\dots+a_k=n}} 
			\frac N{a_k} \int _{\substack{ a_i/N<t_i<(a_i+1)/N \\ 1\le i\le k-1 }} \frac{dt_1\cdots dt_{k-1}}								{t_1\cdots t_{k-1}}\\
		&\gg \frac1{e^{O(k)}}\sum _{\substack{  (a_1,\dots,a_k)\in A^k\\  a_1+\dots+a_k=n}} \frac{N}{a_1\cdots a_k}\ge 					\frac{\alpha_v}{ e^{O(v)}}  \left(  \int_{t\in T}\frac{dt}{t} \right)^k,
}
as desired. Hence we can take $\lambda_3=3\lambda_2$ and $\tau_v\gg \alpha_v e^{-O(v)}$.
\end{proof}


\begin{proof}[Proof that Hypothesis T implies Hypothesis A]
Let $N\ge C_2v^2/\lambda_2$ and $A\subset(N/ev,N/u]$ with $\sum_{a\in A}1/a\ge(1+\lambda_2)/u$. Set $T=\bigcup_{a\in A} (a/N,(a+1)/N)$, so that
\[
\int_T\frac{dt}t=\sum_{a\in A}\frac1a+O\left(\frac1{N/v}\right) 
	= \left(1+O\left(\frac{v^2}N\right)\right)\sum_{a\in A}\frac1a\ge\frac{1+\lambda_2/2}u,
\]
provided that $C_2$ is large enough. If $t_1+t_2+\dots+t_k=1$ and $a_i=[Nt_i]$, then $N-k\leq a_1+\cdots + a_k \le N$. Now, we have that
\[
\frac 1{a_1\dotsm a_k} \geq \idotsint _{\substack{ a_i/N<t_i<(a_i+1)/N \, \forall i }} \frac{dt_1dt_2\dotsm dt_{k-1}}{t_1t_2\cdots t_{k-1}\cdot (Nt_k)},
\]
and so
\[
\sum_{\substack{  (a_1,\dots,a_k)\in A^k\\  N-k\leq a_1+\cdots + a_k \leq N}} \frac 1{a_1\dotsm a_k} \geq
\idotsint\, _{\substack{   t_1+t_2+\dots+t_k=1\\t_1,t_2,\dots, t_k\in T}}\frac{dt_1dt_2\dotsm dt_{k-1}}{N t_1t_2\cdots t_k}\gg   \frac {\tau_v}N \left(  \int_{t\in T}\frac{dt}{t} \right)^k,
\]
provided that $\lambda_2\ge2\lambda_3$. The result follows by averaging over the subsums with $a_1+\cdots + a_k=n$ for each integer $n$ in the interval $N-k \leq n\leq N$. Hence we can take $\lambda_2=2\lambda_3$ and $\alpha_v\gg \tau_v/e^{O(v)}$.
\end{proof}


\begin{proof}[Proof that Hypothesis P implies Hypothesis A] 
Given our set $A$, let $\CP$ be the set of primes in
$\bigcup_{a\in A} (e^a, e^{a+1})$, and let $x=e^{[N]+1}$. Then
\[\sum_{p\in \CP} \frac1p =\sum_{a\in A}\frac1a+O\left(\frac1{N/v}\right)=\left(1+O\left(\frac{v^2}N\right)\right)\sum_{a\in A} \frac1a \ge \frac{1+\lambda_2/2}u,\]
provided that $C_2$ is large enough. So, if we choose $\lambda_2=2\lambda_1$, then we can apply Hypothesis P. Now for each $k$-tuple of primes $(p_1,\dots,p_k)\in \CP^k$ such that $x/2 < p_1\cdots p_k<x$, let $a_j=[\log p_j]$ for all $j$, so that
\[
\sum_{j=1}^k a_j > \sum_{j=1}^k (\log p_j-1) \geq \log (x/2) - k  = [N]+1-\log 2-k
\]
and
\[\sum_{j=1}^k a_j  \leq \sum_{j=1}^k \log p_j  < \log x = [N]+1.
\]
Since $a_j$ are integers, this implies that $N-k \leq \sum_{j=1}^k a_j \leq N$. Hence, noticing that $1/a \sim \sum_{e^a \leq p \leq e^{a+1}} 1/p$, we deduce that
\als{
\sum_{\substack{  (a_1,\dots,a_k)\in A^k\\  N-k\leq  a_1+\dots+a_k\leq N }} \frac 1{a_1\dotsm a_k}  &\gg
e^{-O(k)} \sum_{\substack{  (p_1,\dots,p_k)\in \CP^k\\  x/2 < p_1\cdots p_k<x}} \frac 1{p_1\dotsm p_k} \\
&\gg \frac{e^{-O(k)} \pi_v}{\log x} \cdot \left(\sum_{\substack{ p\in\CP}}\frac1p\right)^k\
\gg\frac{\pi_v}{e^{O(k)}}\cdot\frac{1}{N} \left(  \sum_{\substack{ a\in A  }}  \frac 1a \right)^k .
}
Hence we can take $\lambda_2=2\lambda_1$ and $\alpha_v\gg \pi_v e^{-O(v)}$.
\end{proof}


\begin{proof}[Proof that Hypothesis A  implies Hypothesis P] 
Let $\rho = 1+\delta/(2ev), \eta = \min\{1, \lambda_2\}$ and $N=\log_{\rho} x-ev$. For each integer $a\in[N/ev+1,N/u]$ define
\[
S_a=\sum_{\substack{ p\in \CP\\  \rho^a\le p< \rho^{a+1}}} \frac 1p.
\]
Huxley's prime number theorem for short intervals (see Theorem 10.5 in \cite{ik} and the subsequent discussion) yields $|\{y \leq p \leq y+h\}| \sim h/\log y$ for $y^{7/12+\epsilon} \leq h \leq y$. This implies that
\eq{eq5.1}{
S_a \leq \frac{1+\eta/10}{a} \ll \frac{v}{N} 
	\asymp \frac{\delta}{\log x} 
}
when $C_1$, and thus $x^{1/ev}$, is large enough.

Let $A$ be the set of integers $a$ for which
$S_a\geq  \frac {\eta}{4 a\log (ev)}\sum_{\substack{ p\in \CP  }}  \frac 1p$. Then
\als{
\sum_{\substack{ p\in \CP  \\ [\log_\rho p]\not\in A}} \frac 1p 
	\le O\left(\frac {\delta}{\log x}\right)+\sum_{\substack{ a \not\in A \\ \frac{N}{ev}+1\le a\leq \frac{N}{u} }}  
		\frac {\eta}{4a \log(ev)} \sum_{\substack{ p\in \CP  }}  \frac 1p 
	&\le \left(\frac{\eta}4+O\left(\frac{\delta v}{\log x}\right)\right) \sum_{\substack{ p\in \CP  }}  \frac 1p \\
	&\le\frac {\eta}3 \sum_{\substack{ p\in \CP  }}  \frac 1p,
}
provided that $C_1$ is large enough and $\lambda_1\le100\lambda_2$. Using \eqref{eq5.1} we find that
\[
\sum_{a\in A} \frac1a\ge \left(1-\frac{\eta}{10}\right) \sum_{a\in A}S_a= \left(1-\frac{\eta}{10}\right)\sum_{\substack{ p\in \CP\\ [\log_\rho p]\in A}} \frac 1p \ge\left(1-\frac{\eta}2\right) \sum_{\substack{ p\in \CP  }}  \frac 1p.
\]
So setting $\lambda_1=4\lambda_2$ allows us to apply Hypothesis A. Now for each solution to $a_1+\dotsb+a_k = n\in [N-k, N]$ with
$a_1,\dots,a_k\in A$, consider the primes $p_j\in \CP$ with $[\log_\rho p_j]=a_j$. Note that
$a_j\leq \log_{\rho} p_j<a_j+1$ and so $\log_\rho x-ev-k \leq n\leq \log_\rho (p_1p_2\dotsm p_k) <n+k \leq \log_\rho x$, which implies that $x > p_1 \dotsm p_k > x \rho^{-2ev} = x (1 + \delta/(2ev))^{-2ev} \geq (1-\delta)x$. Hence
\als{
\sum_{\substack{  (p_1,\dots,p_k)\in \CP^k\\  (1-\delta)x<p_1\cdots p_k<x}} \frac 1{p_1\dotsm p_k}
		&\geq \sum_{\substack{  (a_1,\dots,a_k)\in A^k\\  a_1+\dots+a_k=n}} S_{a_1}\cdots S_{a_k}\\
		&\geq \left( \frac {\eta}{4 \log(ev)} \sum_{\substack{ p\in \CP  }} 
		 \frac 1p \right)^k \sum_{\substack{  (a_1,\dots,a_k)	\in A^k\\  a_1+\cdots+a_k=n}} 
		 	\frac 1{a_1\cdots a_k} \\
		&\ge \left( \frac{\eta(1+\lambda_1)}{4u \log(ev)}\right)^k \cdot \frac{\alpha_v} {N} 
			\left(\sum_{\substack{ a\in A  }}  \frac 1a \right)^k  
		\gg \frac{\eta^v\alpha_v}{v^{O(v)}}\cdot\frac{\delta}{\log x} 	
		\left( \sum_{\substack{ p\in \CP  }}  \frac 1p \right)^{k},
}
which proves the desired result.
\end{proof}


\section{Lemmas in additive combinatorics}\label{add-comb}

Let us first introduce some notation. Given two additive sets $A$ and $B$, define the sum set 
$A + B=\{a + b: a \in A, b \in B\}$, the $k$-fold sum set $k A = \{a_1+\dotsb + a_k : a_i \in A\}$, and for any $E\subseteq A\times B$ define the restricted sum set
\[
A{\overset{E}{+}}B=\{a+b:(a,b)\in E\}.
\]
Write also
\[
r_{kA}(n) = \# \{(a_1, \dotsc, a_k) \in A^k: a_1+\dotsb+a_k = n \}
\]
for the number of representations. Finally, a set of the form
\[
P = \{ x_0 + l_1 x_1 + \dotsb + l_d x_d : 0 \leq l_j \leq L_j \text{ for all $j$}\}
\]
is called a {\it generalized arithmetic progression of rank} $d$.

We need three lemmas from additive combinatorics. The first one lets us pass from a restricted sum set to a regular sum set.


\begin{lem}\label{lem6.1}
Let $(G, +)$ be an abelian group. If $E \subseteq A \times A$ satisfies
\[
|E| \ge (1-\delta^2)|A|^2 \quad \text{and} \quad |A{\overset{E}{+}}A| \le K|A|,
\]
then there exists a set $A' \subseteq A$ such that
\[
|A'| \ge (1-\delta)|A| \quad \text{and} \quad |A'-A'| \le \frac{K^2}{1-2\delta}|A|.
\]
\end{lem}


\begin{proof}
This is a variant of the Balog-Szemer\'edi-Gowers theorem (see \cite[Theorem 2.29]{tv}) which can be proved by incorporating the hint for \cite[Exercise 2.5.4]{tv} to the proof of the Balog-Szemer\'edi-Gowers theorem in \cite[Section 6.4]{tv}. We provide a proof for completeness.

Choose
\[
A' = \{a \in A \ \vert \ (a, b) \in E \text{ for at least
  $(1-\delta)|A|$ of $b \in A$}\}.
\]
Now
\[
(1-\delta^2)|A|^2 \leq |E| \leq |A'| |A| + (|A|-|A'|)(1-\delta)|A| \implies |A'| \geq (1-\delta)|A|.
\]

Fix a pair $(a_1, a_2) \in A' \times A'$ and note that $(a_1, b) \in E$ for at least
$(1-\delta)|A|$ of $b \in A$ and $(a_2, b) \in E$ for at least $(1-\delta)|A|$ of $b \in A$. 
Hence there are at least $(1-2\delta)|A|$ elements $b \in A$ for which both $(a_1, b) \in E$ and $(a_2, b) \in E$. Since $a_1-a_2 = (a_1+b) - (a_2+b)$, writing $x = a_1+b$ and $y = a_2+b$, we have
\[
|\{(x, y) \in (A{\overset{E}{+}}A)^2 \ \vert \ x-y = a_1-a_2\}| \geq
(1-2\delta) |A|.
\]
Since the total number of triples $(x, y) \in (A{\overset{E}{+}}A)^2$ is
at most $K^2|A|^{2}$, the claim follows by summing over elements of $A'-A'$.
\end{proof}


The second lemma shows that if $3A$ is small, then we can find a popular large generalized arithmetic progression inside it.


\begin{lem}\label{lem6.2}
Let $K \geq 1$ and let $A$ be a finite subset of the integers such that $|3A| \leq K |A|$. Then there is a generalized arithmetic progression $P \subseteq 3A$ of rank $O_K(1)$ such that $|P| \gg_K |A|$ and $r_{3A}(n) \gg_K |A|^2$ for all $n \in P$. 
\end{lem}


\begin{proof}
This is a variant of the Ruzsa-Chang theorem (see \cite[Theorem 5.30]{tv}). Similarly to that theorem, this can be reduced to the following similar result in $\Z_N := \Z / N\Z$ through theory of Freiman morphisms. 
\end{proof}


\theoremstyle{plain}
\newtheorem*{lem 6.2*}{Lemma \ref{lem6.2}*}

\begin{lem 6.2*}
Let $K, N \geq 1$ with $3 \nmid N$, $\delta > 0$ and let $A$ be a finite subset of $\Z_N$ such that $|3A| \leq K |A|$ and $|A| \geq \delta N$. Then there is a generalized arithmetic progression $P \subseteq 3A$ of rank $O_{K, \delta}(1)$ such that $|P| \gg_{K, \delta} |A|$ and $r_{3A}(n) \gg_{K, \delta} |A|^2$ for all $n \in P$. 
\end{lem 6.2*}


\begin{proof}
This is a consequence of \cite[Theorem 4.43]{tv}, except we have added the requirement $r_{3A}(n) \gg_K |A|^2$ for all $n \in P$ which the proof easily gives. For completeness we sketch the proof.

Let us consider $r_{3A}(n) = 1_A \ast 1_A \ast 1_A(n)$. One has $\sum_{x \in \Z_N} r_{3A}(x) = |A|^3$ whereas $r_{3A}$ is supported on the set $3A$ of cardinality at most $K|A|$. Hence there is $x_0$ such that $r_{3A}(x_0) \geq |A|^2/K$. By translating $A$ by $x_0/3$, we can assume $x_0 = 0$.

Now, writing $\widehat{g}(k) = \sum_{x \in \Z_N} g(x)e(-\tfrac{kx}{N})$ for the Fourier transform,
\als{
|r_{3A}(x)-r_{3A}(0)| &= \frac{1}{N}\left|\sum_{\xi \in \Z_N} \widehat{1_A}(\xi)^3 (e(\tfrac{\xi x}{N})-1)\right| \\
	& \leq \sup_{\xi \in \Z_N} |\widehat{1_A}(\xi)||(e(\tfrac{\xi x}{N})-1)| \cdot \frac{1}{N}\sum_{\xi \in \Z_N} |\widehat{1_A}(\xi)|^2 \leq 2\pi |A| \sup_{\xi \in \Z_N} |\widehat{1_A}(\xi)|\left\Vert\frac{\xi x}{N}\right\Vert,
}
by Parseval's identity and where we write $\Vert y \Vert$ for the distance from the nearest integer.
Hence $r_{3A}(x) \geq |A|^2/(2K)$ for every $x$ in the set
\[
\left \{x \in \Z_N : \ \sup_{\xi \in \Z_N} |\widehat{1_A}(\xi)|\left\Vert\frac{\xi x}{N}\right\Vert < \frac{|A|}{4\pi K}\right\} 
	\supseteq  \left\{x \in \Z_N : \ \sup_{_{\substack{ \xi \in \Z_N \\ |\widehat{1_A}(\xi)| \geq |A|/(4\pi K)  }}} 
		\left\Vert\frac{\xi x}{N}\right\Vert < \frac{1}{4\pi K}\right\} .
\]

By the Fourier concentration lemma \cite[Lemma 4.36]{tv}, there is $d = O_{K, \delta}(1)$ and a set $S = \{\eta_1, \dotsc, \eta_d\} \subset \Z_N$ such that
\[
\{\xi \in \Z_N : \ |\widehat{1_A}(\xi)| \geq |A|/(4\pi K)\} \subseteq \left\{\sum_{j = 1}^d \alpha_j \eta_j : \ \alpha_j \in \{-1, 0, 1\}\right\},
\]
and hence by the triangle inequality $r_{3A}(x) \geq |A|^2/(2 K)$ for every $x$ in the set
\[
\left\{x \in \Z_N : \ \sup_{\xi \in S} \left\Vert\frac{\xi x}{N}\right\Vert < \frac{1}{4 \pi dK}\right\}.
\]
This is a so-called Bohr set which contains the claimed arithmetic progression by \cite[Proposition 4.23]{tv}, the proof of which uses Minskowski's second theorem.
\end{proof}


The third lemma shows that sum sets of generalized arithmetic progressions have large popular subsets.


\begin{lem}\label{lem6.3}
Let $L_j$ be positive integers for $j =1, \dotsc, d$ and let
\[
P = \{ x_0 + l_1 x_1 + \dotsb + l_d x_d : |l_j| \leq L_j \text{ for all $j$}\}
\]
be a generalized arithmetic progression. Let $k\in\mathbb{N}$, $\delta \in (0, 1/6^d), \rho = 1-3\delta^{1/d} \geq 1/2$ and
\[
Q_k = \{kx_0 + l_1 x_1 + \dotsb + l_d x_d : |l_j| \leq \rho k L_j \text{ for all $j$}\}.
\]
Then, for any $n \in Q_k$, one has $r_{kP}(n) \ge (\delta |P|)^{k-1}$.
\end{lem}


\begin{proof}
We proceed by induction on $k$. The case $k = 1$ is trivial, so we assume that the claim holds for some $k\in \mathbb{N}$. Let $n \in Q_{k+1}$, so that $n = (k+1) x_0 +l_1 x_1 + \dotsb + l_d x_d$ with $|l_j| \le \rho (k+1) L_j$ for all $j$. Now
\als{
r_{(k+1)P}(n) &= \sum_{|i_j| \leq L_j \, \forall j} r_{kP}(kx_0 + (l_1-i_1) x_1 + \dotsb + (l_d-i_d) x_d) \\
		& \geq (\delta |P|)^{k-1} \cdot \ \# \{(i_1, \dotsc, i_d):\ |i_j|\leq L_j \ \text{and} \ |l_j-i_j|\leq \rho k L_j 
		\text{ for all $j$} \}.
}
The right hand side is smallest when $l_j=[\rho(k+1)L_j]$ for all $j\in\{1,\dots,d\}$, in which case
\als{
&\#\{i_j\in[-L_j,L_j]:|l_j-i_j|\le\rho k L_j\} \\
&=\#\{i_j:[\rho(k+1)L_j]-\rho k L_j \le i_j \le \min\{L_j, [\rho(k+1)L_j]+\rho k L_j\}\}.
}
Since $2\rho k \geq 2\rho \geq 1$, the minimum above is $L_j$ and hence the number of counted $i_j$ is at least $(1-\rho)L_j$. Hence, since $|P| \leq \prod_{j=1}^d (2L_j+1)\leq \prod_{j=1}^d (3L_j)$, we have that
\[
r_{(k+1)P}(n) \ge (\delta|P|)^{k-1} (1-\rho)^d \prod_{j=1}^d L_j = (\delta|P|)^{k-1} \delta \prod_{j=1}^d (3L_j) \geq (\delta|P|)^k.
\]
\end{proof}


\section{The proof of Hypothesis A}\label{hypa-proof}

Our additive combinatorial tools do not involve logarithmic weights, so instead of Hypothesis A we apply them to prove the following variant.


\begin{thm}\label{thm7.1}
There exists a constant $c_0 > 1$ such that if $1 \leq u \leq v$ and $B$ is a subset of the integers in $(\frac{N}{ev},\frac{N}{u}]$ for which
\[
|B|>\frac{c_0 N}{u^2},
\]
then there is an integer $k\in[u,ev]$ such that
\[
|\{(b_1, \dotsc, b_k) \in B^k : N-k \le b_1 + \dotsb + b_k \le N\}|\gg\frac1{e^{O(k)}k^{k}}\frac{|B|^k}{N}.
\]
\end{thm}


\begin{proof}[Proof that Theorem \ref{thm7.1} implies Hypothesis A with $\lambda_2=2 c_0-1$ and $\alpha_v=v^{-4v}$]
Let $A\subset[N/(ev),N/u]$ as in Hypothesis A. We claim that there must exist $t\in[u,ev]$ such that
\eq{eq7.1}{
\sum_{\substack{ a\in A \\N/(ev)<a\leq N/t}} 1 \ge\frac{c_0N}{t^2}.
}
Indeed, if this is not the case, then
\als{
\sum_{_{\substack{ a\in A \\N/(ev)<a\leq N/u}}}\frac{1}{a} &= -\int_{u}^{ev} \frac{t}{N} d\left(\sum_{_{\substack{ a\in A \\N/(ev)<a\leq N/t }}} 1\right)\\
& = \frac{u}{N} \sum_{_{\substack{ a\in A \\N/(ev)<a\leq N/u }}} 1 + \int_u^{ev}\left(\sum_{_{\substack{ a\in A \\N/(ev)<a\leq N/t }}} 1 \right)\frac{dt}N \\
&< \frac{u}{N} \cdot \frac{c_0 N}{u^2} + c_0 \int_u^{ev}\frac{1}{t^2} dt < 2 \cdot\frac{c_0}{u} = \frac{1+\lambda_2}{u}
}
which is a contradiction. So there is some $t\in[u,ev]$ for which \eqref{eq7.1} holds.

Now, set $B=A\cap (N/(ev),N/t]$, so that the hypothesis of Theorem \ref{thm7.1} is satisfied with $t$ in place of $u$. Let $k$ be as in the conclusion of Theorem \ref{thm7.1}, which necessarily lies in $[u,ev]$. Let $n$ be an integer in $[N-k,N]$ whose number of representations as $b_1+\cdots+b_k$ is maximal. So $n$ has at least $\gg  e^{-O(k)}k^{-k} |B|^k/N$ such representations by Theorem \ref{thm7.1}. Since each $b\in B$ satisfies $1/b\geq u/N \geq 1/N$, and as $|B|\geq c_0 N/t^2\geq c_0N/(ev)^2$, we deduce that
\[
\sum_{\substack{  (b_1,\dots,b_k)\in B^k\\  b_1+\dots+b_k=n}} \frac 1{b_1\dotsm b_k} \gg
\frac 1N \left( \frac{e^{-O(1)} |B|}{Nk}  \right)^k
\gg \frac{e^{-O(v)}}{v^{3v}\cdot N}  \ge \frac{\alpha_v}{N}  \left(  \sum_{\substack{ a\in A  }}  \frac 1a \right)^k
\]
with $\alpha_v \gg v^{-4v}$, since $\sum_{a\in A}\frac1a \ll \log(ev)$.
\end{proof}

\begin{proof}[Proof of Theorem \ref{thm7.1}]
Let $c_0$ be a large positive constant to be determined later. Notice first that if $|B| \leq c_0$, then we only need to find one sum in the interval $[N-k, N]$. In this case the elements of $B$ have size $\leq N/u \leq u$ by the lower bound for $|B|$, so the claim follows trivially.

From now on we can assume that $|B| > c_0$. We claim that if $|B|>c_0N/u^{2}$, then there exists $k$ such that
\[
|\{(b_1, \dots, b_k) \in B^k : N-k \leq b_1 + \dots + b_k \le N\}|\ge c\frac{\delta^{2k}|B|^k}{N k^{k}}
\]
for some appropriate small positive constants $c$ and $\delta$. If $1\le u\le c_0$, our claim is trivial, since there are no sets $B\subset(N/ev,N/u]$ with $|B|>c_0N/u^{2}$.

Now we prove that if the claim holds when $2^{j-1}c_0 \le u \le 2^jc_0$ for some
$j \geq 0$, then it holds when $2^jc_0 \le u \le 2^{j+1}c_0$. Take
\[
E = \left\{(b_1, b_2) \in B \times B :\ r_{2B}(b_1+b_2) \ge \delta^2 \frac{|B|^2}{|2B|}\right\},
\]
so that
\[
|E| \ge |B|^2 - |2B| \cdot \delta^2 \frac{|B|^2}{|2B|} = (1-\delta^2)|B|^2.
\]
Write $C = B {\overset{E}{+}} B \subseteq (2N/(ev), 2N/u]$. We split the rest of the argument into two cases according to whether $|C|>4|B|$ or not.


Consider first the case $|C| > 4|B|$. Then
\[
|C| > \frac{c_0 N}{(u/2)^{2}},
\]
and by induction hypothesis there is an integer $k/2 \in [u/2, ev/2]$ such that
\[
|\{(c_1, \dots, c_{k/2}) \in C^{k/2} : N-k/2 \le c_1 + \dotsb + c_{k/2} \le N\}| \ge c \frac{\delta^k |C|^{k/2}}{N (k/2)^{k/2}}.
\]
Hence by the definition of $C$ we have
\als{
&|\{(b_1, \dots, b_{k}) \in B^{k} : N-k/2 \leq b_1 + \dotsb + b_{k} \leq N\}| \\
		&\ge c\frac{\delta^k |C|^{k/2}}{N (k/2)^{k/2}} 
		\cdot \left(\delta^2 \frac{|B|^2}{|2B|}\right)^{k/2} 
		=  c\frac{\delta^{2k} |B|^{k}}{N {k}^{k/2}} \cdot \left(\frac{2|C|}{|2B|}\right)^{k/2}.
}
The claim now follows if $|C|>|2B|/2k$; this easily follows, since $|C| > 4|B|$ and $|2B| \leq 2N/u \leq 2|B|u/c_0 \leq 2|B| k$.

On the other hand, if $|C| \le 4 |B|$, then by Lemma \ref{lem6.1} there is $B' \subseteq B$ such that $|B'| \geq |B|/2$ and $|B'-B'| \leq 20|B'|$. Then by \cite[Proposition 2.26]{tv} we have $|3B'| \ll |B'|$ and hence Lemma \ref{lem6.2} implies that there is a generalized arithmetic progression $P \subseteq 3B' \subseteq 3B$ of rank $d \ll 1$ such that $|P| \gg |B|$ and $r_{3B}(n) \gg |B|^2$ for all $n \in P$. ``Centralizing'' $P$ we can assume that it is of form
\[
P = \{ x_0 + l_1 x_1 + \dotsb + l_d x_d : |l_j| \leq L_j \text{ for all $j$}\}
\]
for some positive integers $L_j$ (doing this reduces the size of $P$ at most by a factor of $1/3^d \gg 1$). Set $\rho=1-(3\delta)^{1/d}$ and
\[
Q = \{ x_0 + l_1 x_1 + \dotsb + l_d x_d : |l_j| \leq \rho L_j \text{ for all $j$}\} \subset P.
\]
If $\delta<1/(6^d)$, then Lemma \ref{lem6.3} implies that $r_{kP}(n) \geq (\delta |P|)^{k-1}$ for all $k \geq 1$ and $n \in k Q$. Moreover, if $\delta$ is small enough, then $|Q| \geq |P|/3^d \gg |B| > c_0N/u^2$, so that if $c_0$ is large enough, then $|Q| \geq 15 N/u^2$ and, because
\[Q \subseteq P \subseteq 3B \subseteq (3N/(ev), 3N/u],\]
we have
\[
\sum_{q \in Q} \frac 1q \geq \frac{15N}{u^2} \cdot \frac u{3N} \geq \frac 5u.
\]

From here we argue much as in the proof of Proposition \ref{prop5.1}(ii). We begin by removing $[3N/u]$ from $Q$ if it is an element so that $Q\subseteq (3N/(ev), 3N/u-1]$. Let $T=\bigcup_{q\in Q} \left( \frac qN , \frac{q+1}N \right) \subseteq (\frac{3}{ev}, \frac 3u)$, so that
\[
\int_T \frac{dt}t \geq  \sum_{q \in Q} \left( 1 - \frac 1q  \right) \frac 1q  \geq \frac{2}{3} \sum_{q \in Q} \frac 1q  \geq \frac 3u
\]
since $q \geq 3$ for every $q \in Q \subseteq 3B$. Hence, Bleichenbacher's theorem implies that there exists an integer $k/3 \in [u/3, ev/3]$ and $t_1,\ldots, t_{k/3}\in T$ for which $t_1+\ldots + t_{k/3}=1$.  If $t_i\in \left( \frac {q_i}N , \frac{q_i+1}N \right)$ for each $i$ then
$n:=q_1+\ldots + q_{k/3}\in  \frac{k}{3} Q \cap [N-k/3, N]$. Recalling that $r_{3B}(m) \gg |B|^2$ for every $m \in P$ and $r_{\frac{k}{3}P}(\ell) \geq (\delta |P|)^{k/3-1}$ for every $\ell \in \frac{k}{3}Q$, we get that
\als{
&\# \{(b_1, \dotsc, b_{k}) \in B^{k} : N - k/3 \leq b_1 + \dotsb + b_{k} \leq N\} \\
	&\ge \# \{ (b_1, \dotsc, b_{k}) \in B^{k} :  b_1 + \dotsb + b_{k} =n\} \\
	&\geq \sum_{\substack{ m_1+\ldots +m_{k/3} =n \\ m_1,\ldots ,m_{k/3}\in P }}\  \prod_{j=1}^{k/3} r_{3B}(m_j)
		\geq r_{\frac k3\cdot P} (n) \cdot (e^{-O(1)} |B|^2)^{k/3}  \geq
		(\delta |P|)^{k/3-1} \cdot e^{-O(k)} |B|^{2k/3} \\
	&\ge e^{-O(k)}(\delta |B|)^{k-1} \ge ue^{-O(k)}(\delta |B|)^k/N,
}
as $P\subseteq 3B$.
The result follows, since the right hand side is $\gg e^{-O(k)} |B|^k/N$ for every fixed $\delta > 0$.
\end{proof}

\begin{rem}\label{rem7.1}
One could compute the constant $c_0$ explicitly and thereby the constants $\lambda_i$ in Hypotheses P, A and T and, eventually, $\lambda$ in Theorem \ref{thm1}. However, $c_0$ will be relatively large, since the implied constants in Lemma \ref{lem6.2} are rather large. If one is interested in optimizing $\lambda$, one could, instead of Lemma \ref{lem6.2}, use a result of Lev \cite{lev} to show that if the number of ``popular'' elements in $|2A|$ is at most $K|A|$ for some $K < (3+\sqrt{5})/2$, then $2A$ contains a ``popular'' arithmetic progression (of rank $1$). Modifying the above arguments, this would lead to Theorem \ref{thm1} with a smaller and more easily calculable $\lambda < 21$. However, this argument would not yield Hypotheses P, A and T when $u$ is not close to $1$ and, in particular, not the latter conclusion in Remark \ref{rem1.4}. By applying Bleichenbacher's theorem in a different way one could probably improve $\lambda$ further, but not to an arbitrarily small constant, as desired.
\end{rem}


\section{Some combinatorial lemmas}\label{comb-lemmas}

We devote this section to proving some combinatorial lemmas we will need in next section where we investigate some further consequences of Hypotheses A, P and T.


\begin{lem}\label{lem3.1}
Let $B$ be a finite subset of the numbers in $(y,z]$ and associate to each $b\in B$ a positive weight $w(b)$.  For any $x\geq z$ there exists a positive integer $k\leq x/y$ such that
\[
\sum_{\substack{ b_1,\dots ,b_k\in B \\b_1+\dots+b_k\in (x-z,x] }} w(b_1)w(b_2)\dotsm w(b_k) \geq \frac 1{x/y}
\left(\sum_{b\in B}w(b)\right)^k  .
\]

In particular, letting $w(b) = 1$ for all $b \in B$ yields: Let $B$ be a finite subset of the numbers in $(y,z]$.  For any $x\geq z$ there exists a positive integer $k$ such that the number of $k$-tuples $b_1,\dots ,b_k\in B$ for which
$b_1+\dots +b_k\in (x-z, x]$ is  $\geq |B|^k/ (x/y)$.
\end{lem}


\begin{proof}
Note that if $b_1+\dots+b_k\in(x-z,x]$ for some $b_1,\dots,b_k\in B$, then $ky<b_1+\dots +b_k\leq x$ and so $1\leq k<x/y$. Therefore if we let $K=[x/y]$, then
\als{
\beta_k&:=\sum_{\substack{ b_1,\dots ,b_k\in B \\  b_1+\dots +b_k\in (x-z, x] }} w(b_1)w(b_2)\dotsm w(b_k) \bigg/
\left(  \sum_{b\in B}  w(b) \right)^k \\
&=\sum_{\substack{ b_1,\dots ,b_K\in B \\  b_1+\dots +b_k\in (x-z, x]}} w(b_1)w(b_2)\dotsm w(b_K) \bigg/
\left(\sum_{b\in B}w(b)\right)^K.
}
Consequently, we find that
\[
\sum_{k=1}^K \beta_k = \left(\sum_{b\in B}w(b)\right)^{-K}\sum_{\substack{ b_1,\dots ,b_K\in B}} w(b_1)w(b_2)\dotsm w(b_K) \sum_{\substack{ 1\leq k\leq K \\ b_1+\dots +b_k\in (x-z,x]}} 1 \geq 1,
\]
since the differences in each sequence $b_1, b_1+b_2, \dots , b_1+\dots +b_K$ are $\leq z$, whereas $b_1\leq z\leq x$ and $b_1+\dots +b_K > yK > x-y > x-z$. Taking the maximum of the $\beta_k$ then yields the desired result.\qed
\end{proof}


\begin{cor}\label{cor3.2}
Let $\CP$ be a  subset of the primes in $(x^{1/ev},x^{1/u}]$ for some $1\le u\le v\le(\log x)/e$. For any $X\geq x^{1/u}$ there exists a positive integer $\ell\leq K:=\frac{ev \log X}{\log x}$ such that
\[
\sum_{\substack{ q_1,\dots ,q_\ell\in \CP \\  Xx^{-1/u}< q_1q_2\dotsm q_\ell\leq X }} \frac{1}{q_1q_2\dotsm q_\ell}  \geq \frac 1{K}
\left(  \sum_{q\in \CP} \frac 1q \right)^\ell.
\]
\end{cor}


\begin{proof}
Apply Lemma \ref{lem3.1} with $B=\{\log q:\ q\in \CP\}\subset ((1/ev)\log x, (1/u)\log x]$ and $w(b)=e^{-b}=1/q$, and then take $\ell = k$ as obtained in that lemma.
\end{proof}


\section{Further remarks on Hypotheses P, A and T}\label{a-p-t}

We conclude our paper with an investigation of some other consequences of Hypotheses P,A and T. 


\begin{prop}\label{prop8.1}
Suppose Hypothesis P holds with $\pi_v=1/v^{O(v)}$. If $\epsilon>0$, $v^2\le \lambda_1 \log x/C_1$, $v \geq 1$ and $\CP$ is a subset of the primes in $(x^{1/ev},x^{1/v}]$ for which $\sigma:=\sum_{\substack{ p\in\CP}}\frac1p \geq \frac {1+\lambda_1+\epsilon}v$, then there exists an integer $k\in[v,ev]$ such that
\[
\Psi(x;\CP) \gg_\epsilon \frac x{\log x} \cdot  \frac{\sigma^{k+O(1/\sigma)}}{k\cdot k!}.
\]
\end{prop}


\begin{proof}
We can clearly assume that $\epsilon$ is small. Also, since $\Psi(x;\CP)\ge1$, we may assume that $x$ is large. If $v<2$ or $\frac{1+\lambda_1}{v-1}>\sigma\geq \frac{1+\lambda_1+\epsilon}{v}$, then $v<\max\{2, \frac{1+\lambda_1+\epsilon}{\epsilon}\} = \frac{1+\lambda_1+\epsilon}{\epsilon}$ and the proposition follows by Hypothesis P. So we may impose the additional assumptions that $v\ge2$ and $\sigma\ge\frac{1+\lambda_1}{v-1}$.

Select $w=(1+\lambda_1)/\sigma$ and note that $w\in [1,v-1]$, since $\sigma\le\sum_{x^{1/(ev)}<p\le x^{1/v}}1/p\lesssim1$. We begin by applying Corollary \ref{cor3.2} with $X=x^{1-\frac wv}\ge x^{1/v}$. So there is $\ell\ll v$ such that
\[
\sum_{\substack{ q_1,\dots ,q_\ell\in \CP \\x^{1-\frac{w+1}v}< q_1q_2\dotsm q_\ell\leq x^{1-\frac wv}}} \frac{1}{q_1q_2\dotsm q_\ell}\gg \frac{\sigma^\ell}v.
\]
For each $q_1\cdots q_\ell$ in the above sum we apply Hypothesis P with $x$ replaced by $x/(q_1\cdots q_\ell)$, and both $u$ and $v$ replaced by $V=\log(x/(q_1\cdots q_\ell))/\log(x^{1/v})$, which is possible since $V\in[w,w+1]$. Consequently,
\[
\sum_{\substack{ (q_1,\dots, q_\ell,  p_1,\dots,p_m )\in \CP^{\ell+m}\\  x/2 < q_1\cdots q_\ell p_1\cdots p_m<x}} \frac 1{q_1\cdots q_\ell p_1\cdots p_m}
\gg \frac{\sigma^{\ell+m}}{w^{O(w)}v\log x}.
\]
Finally, note that
\als{
\Psi(x;\CP)&\ge\frac 1{(\ell+m)!} \sum _{\substack{(q_1,\dots,q_\ell,p_1,\dots,p_m)\in\CP^{\ell+m} \\ q_1\cdots q_\ell p_1\cdots p_m\le x}} 1\\
&\ge\frac x{2(\ell+m)!}\sum_{\substack{ (q_1,\ldots, q_\ell,  p_1,\dots,p_m )\in \CP^{\ell+m}\\  x/2 < q_1\cdots q_\ell p_1\cdots p_m<x}} \frac 1{q_1\cdots q_\ell p_1\cdots p_m}\\
&\gg\frac x{\log x}\cdot \frac{\sigma^{\ell+m}}{w^{O(w)}v(\ell+m)!}.
}
Letting $k=\ell+m$ and observing that, necessarily, $k \in [v, ev]$ and $w \log w \ll \sigma^{-1}\log(1/\sigma)$ completes the proof of the proposition.\qed
\end{proof}

\begin{cor}\label{cor8.2}
There exist constants $c>1$ and $c'>0$ such that if $1\le v \leq c' \sqrt{\log x}$ and $\CP$ is a subset of the primes in $(x^{1/ev},x^{1/v}]$ for which
\[
\sum_{\substack{ p\in \CP  }}  \frac 1p \geq \frac {\max\{c,\log v\}}{v},
\]
then there exists an integer $k\in[v,ev]$ such that
\[
\Psi(x;\CP) \ge \Psi_k(x;\CP) \ge\frac1{k!}\left(e^{-O(1)}\sum_{\substack{ p\in\CP}}\frac1p\right)^k\,\cdot\,\frac x{\log x},
\]
as $x\to\infty$, where $\Psi_k(x;\CP)$ denotes the number of integers $n\leq x$ such that $n$ is squarefree, $n$ has exactly $k$ prime factors, and all of the prime factors of $n$ come from $\CP$.
\end{cor}


We conjecture that Corollary \ref{cor8.2} holds under the weaker assumption that $\sum_{\substack{ p\in\CP}}\frac1p \ge c/v$ for any $c>1$, with the implied constant depending at most on $c$.

Consider more generally $\Psi_k(x, \CP)$ for any $\CP \subset \{p \leq x\}$. If $n\in(\sqrt{x},x]$ is counted by $\Psi_k(x;\CP)$, then we can uniquely write $n=mp$ with $p$ prime and $m$ composed of primes $< p$. Note that $p^k\ge n\ge\sqrt{x}$. So $m$ has $k-1$ prime factors and it is $\le x^{1-1/(2k)}$. Now for each $m$ the number of such primes $p$ is $\lesssim (x/m)/\log (x/m)$ as $x\to\infty$, and therefore
\[
\Psi_k(x;\CP)  \lesssim \sqrt{x}+\frac {2kx}{\log x} \sum_{\substack{ p|m \implies p\in \CP \\ \omega(m)=k-1}} \frac{\mu^2(m)}{m} \leq \sqrt{x}+ \frac{2k}{(k-1)!}\left( \sum_{\substack{ p\in\CP}}\frac1p\right)^{k-1}\cdot\frac x{\log x},
\]
as $x\to\infty$. In particular one cannot significantly improve the lower bound in Corollary \ref{cor8.2}.

It is not difficult to prove corollaries of Hypotheses A and T that are analogous to Proposition \ref{prop8.1}. Thus we have


\begin{prop}\label{prop8.3}  
Suppose that Hypothesis A holds for some $\lambda_2>0$ and $C_2>1$ with $\alpha_v=1/v^{O(v)}$.  If $\epsilon>0$, $v^2\leq\lambda_2 N/C_2$, $v \geq 1$ and $A$ is a subset of the integers in $(\frac{N}{ev}, \frac{N}{v}]$ such that $\alpha:=\sum_{a\in A} 1/a \geq(1+\lambda_2+\epsilon)/v$, then there exists an integer $k\in[v,ev]$, and an integer $n$ in the range $N-k\leq n\leq N$,  such that
\[
\sum_{\substack{  (a_1,\dots,a_k)\in A^k\\  a_1+\dots+a_k=n}} \frac 1{a_1\cdots a_k}
\gg_\epsilon\frac{1}{N} \cdot \frac{\alpha^{k+O(1/\alpha)}}k.
\]
\end{prop}


Similarly, we have the following result.


\begin{prop}\label{prop8.4}  
Suppose that Hypothesis T holds for some $\lambda_3>0$ with $\tau_v=1/v^{O(v)}$.  If $\epsilon>0$, $v \geq 1$ and $T$ is an open subset of $(\frac{1}{ev}, \frac{1}{v}]$ for which $\tau:=\int_{t\in T} dt/t \geq(1+\lambda_3+\epsilon)/v$, then there exists an integer $k\in[v,ev]$ such that
\[
{\int \dots \int}_{\substack{ t_1+t_2+\dots+t_k=1\\t_1,t_2,\dots, t_k\in T}}\frac{dt_1dt_2\dotsm dt_{k-1}}{t_1t_2\cdots t_k}\gg_\epsilon \frac{\tau^{k+O(1/\tau)}}k.
\]
\end{prop}


In the proofs of Propositions \ref{prop8.3} and \ref{prop8.4} we need appropriate analogues to Corollary \ref{cor3.2}. The needed result for the proof of Proposition \ref{prop8.3} follows in a straightforward way from Lemma \ref{lem3.1}. For the proof of Proposition  \ref{prop8.4} we make note, without proof, of the appropriate result:


\begin{lem}\label{lem8.5}  Let $T$ be an open  subset of  $(\frac1{ev},\frac1{u}]$.  For any $w\geq  {1/v}$ there exists a positive integer $\ell\leq evw$ such that
\[
{\int \dots \int}_{\substack{ w-1/v<t_1+t_2+\dots+t_\ell\leq w\\t_1,t_2,\dots, t_\ell\in T}}\frac{dt_1dt_2\dotsm dt_{\ell}}{t_1t_2\cdots t_\ell}\geq \frac 1{evw}
\left( \int_{t\in T} \frac{dt}t  \right)^\ell  .
\]
\end{lem}



\end{document}